\documentclass[11pt]{amsart}
\usepackage[a4paper,left=2cm,right=2cm,bottom=2.5cm]{geometry}
\usepackage[utf8]{inputenc}
\usepackage[english]{babel}
\usepackage{amsmath,amssymb}
\usepackage{amsthm}
 \usepackage{stmaryrd} 
\usepackage{tikz}
\usepackage{tkz-tab}
\usepackage{listings}
\usepackage{hyperref}
\usepackage{cleveref}
\usepackage{color} 
\definecolor{mygreen}{RGB}{28,172,0} 
\definecolor{mylilas}{RGB}{170,55,241}

\newtheorem{theorem}{Theorem}[section]

\newtheorem{proposition}[theorem]{Proposition}

\newtheorem{remark}[theorem]{Remark}

\newcommand{\R}{\mathbb{R}}

\newcommand{\N}{\mathbb{N}}

\newcommand{\D}{\mathcal{D}}
\newcommand{\Om}{\Omega}
\newcommand{\dOm}{{\partial \Omega}}

\newcommand{\eps}{\varepsilon}
\newcommand{\Vol}{\text{Vol}}
\newcommand{\Per}{\text{Per}}
\newcommand{\W}{\text{W}}

\newcommand{\Jac}{\text{Jac}}

\begin{document}
\title{Numerical exploration of the range of shape functionals using neural networks}
\author[]{Eloi Martinet}
\address[Eloi Martinet]{Institute of Mathematics, University of Würzburg, Germany}
\email{eloi.martinet@uni-wuerzburg.de}

\author[]{Ilias Ftouhi}
\address[Ilias Ftouhi]{Laboratoire MIPA, N\^imes University, Site des Carmes, Place Gabriel P\'eri, 30000 N\^imes, France}
\email{ilias.ftouhi@unimes.fr}

\lstset{language=Matlab,%
    breaklines=true,%
    morekeywords={matlab2tikz},
    keywordstyle=\color{blue},%
    morekeywords=[2]{1}, keywordstyle=[2]{\color{black}},
    identifierstyle=\color{black},%
    stringstyle=\color{mylilas},
    commentstyle=\color{mygreen},%
    showstringspaces=false,
    numbers=left,%
    numberstyle={\tiny \color{black}},
    numbersep=9pt, 
    emph=[1]{for,end,break},emphstyle=[1]\color{red}, 
}

\date{\today}

\begin{abstract}
We introduce a novel numerical framework for the exploration of Blaschke--Santaló diagrams, which are efficient tools characterizing the possible inequalities relating some given shape functionals. We introduce a parametrization of convex bodies in arbitrary dimensions using a specific invertible neural network architecture based on gauge functions, allowing an intrinsic conservation of the convexity of the sets during the shape optimization process. To achieve a uniform sampling inside the diagram, and thus a satisfying description of it, we introduce an interacting particle system that minimizes a Riesz energy functional via automatic differentiation in PyTorch. The effectiveness of the method is demonstrated on several diagrams involving both geometric and PDE-type functionals for convex bodies of $\R^2$ and $\R^3$, namely, the volume, the perimeter, the moment of inertia, the torsional rigidity, the Willmore energy, and the first two Neumann eigenvalues of the Laplacian. 
\end{abstract}



\maketitle


\section{Introduction}

Visual representations of the range of vector-valued shape functionals—commonly referred to as Blaschke--Santaló diagrams—offer a geometric framework for understanding all possible inequalities that relate prescribed shape functionals within a given class of sets. The terminology honors the pioneering contributions of W. Blaschke \cite{blaschke} and L. Santaló \cite{santalo}, who initiated the systematic study of inequalities linking three geometric functionals on convex bodies. Since then, Blaschke--Santaló diagrams have been widely investigated, particularly in the setting of convex geometry, and also been extended to encompass triplets combining geometric and spectral quantities. For various examples on purely geometric functionals, we refer to the following non-exhaustive list of works \cite{BHT,cifre3,MR3653891,branden,newnew,delyon2,Gastaldello2023Jul}, and for examples involving both geometric and spectral quantities, we refer to \cite{Antunes2010,MR3068840,BBF,zbMATH07369278,ftouh,ftouhi_henrot,LZ,zbMATH06736468}. Also we note that special attention has been given to diagrams involving the Cheeger constant which corresponds to the first eigenvalue of the $1$-Laplacian operator with Dirichlet boundary conditions, see \cite{ftouhi_inequality,ftouhi_cheeger,zbMATH07837343}.

All the aforementioned works focus on establishing theoretical results concerning these diagrams. This includes deriving explicit characterizations in the case of purely geometric functionals, as well as proving qualitative properties when spectral or PDE-related quantities are involved, situations in which obtaining an explicit description is often highly challenging, if not altogether out of reach. 

To the best of our knowledge, the numerical investigation of these diagrams has received relatively little attention in the literature. A first natural numerical approach was proposed in \cite{freitas_antunes}, where the authors examine the diagram involving the perimeter, the area, and the first Dirichlet eigenvalue by approximating it through point clouds generated from the sampling of thousands of random convex polygons. Subsequently, in \cite{Ftouhi2025Jun}, I. Ftouhi proposed a hybrid approach that combines theoretical results on the vertical convexity of the diagrams with numerical shape optimization techniques to determine their boundaries. This allowed him to consider diagrams involving the first Dirichlet eigenvalue and geometric quantities such as the perimeter, the area, and the diameter of planar convex sets. Unfortunately, this approach relies on a priori knowledge of the vertical convexity of the diagram, a property that may be highly nontrivial to establish. As a consequence, the method is not readily applicable to more general classes of diagrams. Then came a paper \cite{oudet_BS} by B. Bogosel, G. Buttazzo and E. Oudet who proposed a method based on Voronoi tesselation and the use of Llyod's algorithm to "uniformly" cover Blaschke--Santal\'o diagrams. The method was successfully applied to diagrams on matrices involving the trace and determinant functionals as well as to a purely geometric diagram incorporating the moment of inertia, the perimeter, and the area of convex sets with two axes of symmetry. However, as far as we know, it could not be extended to settings involving spectral functionals or more complicated geometries such as general convex sets.

In the present paper, we propose a method based on a new neural network architecture to parametrize convex sets. This parametrization has both the advantages of being smooth and without any constraint. We then show how this new parametrization can be used to uniformly sample Blaschke--Santal\'o diagrams using a system of interacting particles. Finally, we show how our new method performs on different types of Blaschke--Santal\'o diagrams.

\section{Shape representation by an invertible neural network}

\subsection{Previous works}

Shape optimization among convex sets has received a lot of attention in the theoretical shape optimization community. Indeed, it is both a natural set of shapes to consider, and at the same time a class that is "rigid" enough so that proving optimality properties is significantly easier than in more generic classes of shapes. On the other hand, numerical shape optimization needs more care to efficiently and precisely take into account the additional convexity constraint. 

One of the first approaches was given by E. Oudet in \cite{Oudet2004}, where the convexity constraint is weakly enforced through the introduction of a penalization term measuring the discrepancy between the current set and its convex hull. In subsequent work \cite{Lachand-Robert2006Jul}, the same author, together with T. Lachand-Robert suggested representing a convex set as the intersection of half-spaces. By construction, however, this representation precludes the approximation of smooth convex sets. Later, E. Oudet \cite{Oudet2013Mar} proposed representing a convex set through a discretization of its support function. This approach has since been widely adopted; see, for instance, \cite{Antunes2022May,Bayen2012Dec,Bogosel2023Feb2,Bogosel2024Nov,Ftouhi2025Jun}. At last, we emphasize that these techniques are also implemented using gauge functions, which, like support functions, furnish a parametrization of convex sets.

However, although these latter approaches yield satisfactory results, for instance, in the parametrization of sets of constant width, it presents several technical difficulties that make its implementation somewhat challenging for the purposes considered in the present work. The first difficulty stems from the fact that the support (or gauge) function must satisfy a second-order differential inequality in order to characterize a convex set, thereby leading to a constrained optimization problem. Depending on the chosen parametrization, the convexity inequality can be enforced in an exact (i.e., non-relaxed) form \cite{Bayen2012Dec,Bogosel2023Feb2}, thereby ensuring that the admissible shapes are genuinely convex. However, this has so far only been achieved for two-dimensional convex sets, and, when using a Fourier basis discretization, it leads to non-linear constraints on the coefficients. In \cite{Antunes2022May}, the method based on the truncation of the Fourier expansion is generalized to the three-dimensional setting. However, the semi-definite constraint is relaxed to a linear one, which can result in a loss of convexity of the admissible shapes, sometimes to a considerable extent. 

In contrast to previous methods, the approach in \cite{Bartels2020Apr} utilizes a direct triangulation-based shape variation framework, wherein the deformation field is constrained to satisfy a convexity-preserving criterion.

Finally, several neural-network-based representations of convex sets have already been proposed, either by enforcing the convexity of the learned level-set function \cite{Martinet2025May} or by adopting a phase-field formulation \cite{Deng2019Sep}.

\subsection{Sublinear architecture}

In this section, we show how sublinear neural networks allows for a flexible and robust parametrization of convex sets. While this method has been proposed in a previous paper by the first author \cite{martinet2026}, we recall here the main principles and results for the sake of completeness. We refer to the previously cited paper for the proofs of the different statements.

The main idea of the method is to represent convex sets as the image of a bijective map making use of a sublinear function, which is given by a new neural network architecture. Note that the idea of representing a shape as a neural diffeomorphism has already been considered in \cite{BelieresFrendo2025Apr}, where authors utilize symplectic maps to represent geometry, thereby ensuring intrinsic measure conservation.

In what follows, $\mathcal{K}_d$ will denote the class of convex bodies of $\mathbb{R}^d$ that contain the origin in their interior. As the shape functionals considered hereafter are translation-invariant, this assumption entails no loss of generality.

Let $p : \R^d \to \R$ be a positive, continuous, and sublinear function, i.e., $p(\lambda x) = \lambda p(x)$ and $p(x+y) \leq p(x) + p(y)$, for every $\lambda >0$ and $x,y\in \mathbb{R}^d$. Then, the set
$$
    \left\{ x \in \R^d : p(x) < 1 \right\}
$$
is an open and convex domain. Conversely, every open and convex body of $\R^d$ is, up to a translation, of the form above for a certain continuous, positive, and sublinear function which corresponds to its gauge function. Using these properties, we can show that convex sets can be described as homeomorphisms of the ball.
\begin{proposition}  
    Let $p$ be a positive, continuous, and sublinear function. The function $\phi:x\longmapsto \frac{\|x\|}{p(x)}x$ is a homeomorphism from $B$ to $\phi(B)$ and the image set $\phi(B)$ is convex. 
\end{proposition}

In what follows, convex sets will be represented via such homeomorphisms. Now, it is important to be able to faithfully represent any positive, continuous, and sublinear function $p$. To this end, we will derive an architecture of a neural network $p_\theta : \R^d \to \R$  (where $\theta$ represents the parameters) which intrinsically satisfies those constraints. A reasonable choice of architecture would be
$$
    p_\theta(x) = \max_{1 \leq i \leq N} (w_i \cdot x),
$$
with $w_i \in \R^d$ for all $i$ (i.e., a single layer neural network with MaxOut activation function \cite{Goodfellow2013Feb}). This is motivated by the fact that any sublinear function is a supremum of linear forms. Then we can set $\phi_\theta(x) := \frac{\|x\|}{p_\theta(x)}x$ to get that $\phi_\theta(B)$ is convex. However, it may be more convenient to be able to define a neural network that is smooth (at least at the boundary of $B$) in order to compute usual quantities of interest (normals, curvature, etc). A natural idea, used, for example, in \cite{Deng2019Sep}, is to replace the previous representation $p_\theta(x) = \max_{1 \leq i \leq N} (w_i \cdot x)$ by a smoothed version obtained via the LogSumExp (LSE) function:
\begin{equation}
    \label{eq:p_LSE}  
    p_\theta(x) = \beta \text{LSE}(W^T x),
\end{equation}
where $\beta > 0$, $W = (w_1 \cdots w_N) \in \R^{d\times N}$, and
$$
    \text{LSE}(y) := \log \sum_{i=1}^N e^{y_i}.
$$

It is well known that $\alpha^{-1} \text{LSE}(\alpha y)$ converges to $\max_{1 \leq i \leq N} y_i$ as $\alpha \to \infty$. The function $p_\theta$ can be seen as a single-layer neural network with learnable parameters $\theta = \{\beta, W\}$.

Note that in \cref{eq:p_LSE}, the function $p_\theta$ is no longer sublinear. However, we can define a function which agrees with the log-sum-exp on the sphere and is sublinear as one sees in the following proposition. 
\begin{proposition}
    The function
    \begin{equation}
        \label{eq:lse_net}
        p_\theta(x) := \beta \|x\| \text{LSE}\left(W^T \frac{x}{\|x\|}\right)
    \end{equation}
    is continuous, sublinear, and is $C^\infty$ on $\R^n \setminus \{0\}$.
\end{proposition}
The proof of this proposition (and the following one) is provided in \cite{martinet2026}.

The set $\Om_\theta := \phi_\theta(B)$ is convex by the construction of the neural network $\phi_\theta(x) := \frac{\|x\|}{p_\theta(x)}x$. However, it is of core interest to know whether every convex set can be represented using a neural network. Since any $K_\theta$ must be smooth, this would not be the case; however, in \cite{martinet2026}, E. Martinet proves the following ``Universal Approximation Theorem'', which ensures that this approach yields a consistent approximation of convex sets.
\begin{theorem}
    Define 
    $$
        \mathcal{K}_d^\text{NN} := \left\{ \Om_\theta : \beta > 0, W \in \R^{d\times N}, N \in \N \right\}
    $$
    as the set of convex bodies that can be represented by the architecture given by \cref{eq:lse_net}. Then $\mathcal{K}_d^\text{NN}$ is dense in $\mathcal{K}_d$ with respect to the Hausdorff distance.
\end{theorem}

\subsection{Encoding symmetries}

Oftentimes, one may be interested in convex bodies that satisfy certain symmetries. Let $G$ be a finite subgroup of isometries of $\R^n$ and $\rho_\theta:\R^n \to \R$ be some neural network given in the previous section, and define
\begin{equation}
    \label{eq:symmetric}
    \phi_\theta^G(x) := x \frac{\|x\|}{\sum_{g \in G} \rho_\theta(g^{-1}x)}.
\end{equation}
Then, $\phi_\theta^G$ is $G$-equivariant, which means that $g.\phi^G_\theta(x) = \phi_\theta^G(g.x)$. This, in turn, implies that the set $\Om_\theta := \phi^G_\theta(B)$ is $G$-invariant, meaning that $g.\Om_\theta = \Om_\theta$ for all $g \in G$.

\section{Sampling the diagram}

Let us denote by $\D := \left\{ F(\Om) : \Om \in \mathcal{K}_d \right\}$ the Blaschke--Santaló diagram of interest, where $F: \mathcal{K}_d \to \R^n$ is some vector-valued shape functional. In order to uniformly sample the diagram $\D$, we will look for a distribution of repulsive electric charges that are at equilibrium inside the diagram $\D$. Formally speaking, this amounts to minimizing the so-called \textit{Riesz energy}:
\begin{equation}
    \label{eq:riesz}
    \min_{y_1, \dots , y_N \in \D} \sum_{1 \leq i\neq j \leq N} \left|y_i - y_j\right|^{-s} 
    \quad \iff \quad
    \min_{\Om_1, \dots , \Om_N \in \mathcal{K}_d} \sum_{1 \leq i\neq j \leq N} \left|F(\Om_i) - F(\Om_j)\right|^{-s},
\end{equation}
with $s>0$. Note that when $s \to \infty$, \cref{eq:riesz} reduces to the famous \textit{sphere-packing problem}
$$
    \min_{y_1, \dots , y_N \in \D} \max_{1 \leq i\neq j \leq N} \left|y_i - y_j\right|.
$$
With that in mind, it is reasonable to expect that for $s$ large enough, the points will distribute uniformly in $\D$, which turns out to be indeed the case, as shown by \textit{Poppy-seed Bagel Theorem} \cite{borodachov_discrete_2019}. For simplification purposes, we state the theorem in the context of smooth manifolds:
\begin{theorem}[Th. 8.5.2 in \cite{borodachov_discrete_2019}]
    \label{th:poppy-seed}
    Let $A \subset \R^p$ be a smooth, compact submanifold of dimension $d<s$. Let $\omega_N = (y_1,\dots, y_N) \in A^N$ be a minimizing configuration of 
    $$
        \min_{(y_1,\dots, y_N) \in A^N} \sum_{1 \leq i\neq j \leq N} \left|y_i - y_j\right|^{-s}.
    $$
    Then, the measure $\nu(\omega_N) := \frac{1}{N} \sum_{y\in \omega_N} \delta_y$ weakly converges to the uniform probability measure on $A$.
\end{theorem}
In the following, we formulate a refined version of the theorem, ensuring that the minimal pairwise distance between points is bounded below by a positive constant. 
\begin{theorem}[Th. 8.8.1 in \cite{borodachov_discrete_2019}]
    Under the assumptions of the previous theorem and the fact that $A$ is path-connected, there exists a constant $C>0$ such that for all $N \geq 2$,
    $$
        \min_{\substack{x,y \in \omega_N \\ x \neq y}} |x-y| \geq \frac{C}{N^{1/d}}.
    $$
\end{theorem}
Using those theorems as motivation for our method (although, we do not necessarily satisfy the assumptions), we will numerically approximate \cref{eq:riesz} by
\begin{equation}
    \label{eq:main-opt}
    \min_{\theta_1, \dots , \theta_N} \sum_{1 \leq i\neq j \leq N} \left|F(\Om_{\theta_i}) - F(\Om_{\theta_j})\right|^{-s}.
\end{equation}

Using PyTorch \cite{Paszke2019Dec}, this problem can be easily solved by some l-BFGS method, as long as the derivative of the function $\theta \mapsto F(\Om_\theta)$ can be automatically computed. We shall assume this condition to hold throughout the remainder of the paper.

\section{Numerical computation of the involved shape functionals}

In this section, we describe in detail the computation of several common shape functionals, highlighting the advantages offered by PyTorch’s automatic differentiation capabilities.

\subsection{Integral quantities}\label{ss:integral_quantities}

For integral quantities, the idea is simply to formulate the computations on the reference domain. To this purpose, we use the volumic and surfacic changes of variable formulas 
$$
    \int_{\Om_\theta} f dx
    = \int_B (f \circ \phi_\theta) \Jac(\phi_\theta) dx
    \quad \mbox{ and } \quad
    \int_{\partial \Om_\theta} g d\sigma
    = \int_{\partial B}(g \circ \phi_\theta) \Jac_{\partial B}(\phi_\theta)d\sigma,
$$
where 
$$
    \Jac(\phi_\theta) = |\det(D\phi_\theta)|
    \quad \mbox{ and } \quad
    \Jac_{\partial B}(\phi_\theta) = \Jac(\phi_\theta)|\left(D\phi_\theta\right)^{-T} n_B|,
$$
and $n_B(x) := x$ on $\partial B$ is the unit outward normal vector. We will apply those formulas to compute the volume $\Vol$, perimeter $\Per$ and moment of inertia $W$ of a shape $\Om_\theta$.

\begin{remark}
    Note that nowhere we actually need the convexity or special form of $\phi_\theta$. This approach can and will be extended to an arbitrary invertible neural network $\phi_\theta$ is a forthcoming paper.
\end{remark}

For the volume, we can simply compute
$$
    \Vol(\Om_\theta)= \int_B \Jac \phi_\theta dx.
$$

We proceed similarly for the moment of inertia, defined as
$$
    W(\Om_\theta) := \int_\Om |x - x_{\Om_\theta}|^2 dx,
$$
where $x_{\Om_\theta} := \int_{\Om_\theta} x dx$ is the center of mass of $\Om_\theta$. By a change of variables, the latter can be expressed as
$$
    x_{\Om_\theta} =  \int_B \phi_\theta(x) \Jac \phi_\theta(x) dx,
$$
while the former is
$$
    W(\Om_\theta) = \int_B |\phi_\theta(x) - x_{\Om_\theta}|^2\Jac \phi_\theta(x) dx.
$$

In the same fashion, the perimeter can be pulled back to the boundary of the reference domain and is computed as follows
$$
    \Per(\Om_\theta) = \int_{\partial B} \Jac_{\partial B}(\phi_\theta)d\sigma.
$$

\subsection{Curvature-related quantities}

The mean curvature of a smooth set $\Om$ at $y \in \dOm$ can be defined as 
$$
    H_\Om(y) := \frac{1}{d-1} \text{div} n_\Om (y) = \frac{1}{d-1} \text{Tr}(D n_\Om(y)),
$$
where $n_\Om$ is a unit length extension of the normal vector to $\R^d$, see for example \cite[Proposition 5.4.8]{HPb}. This quantity appears in the \textit{Willmore energy}:
$$
    E(\Om) := \int_\dOm H_\Om^2(y) dy.
$$
In the case of a parameterized set $\Om_\theta$, the normal vector $n_{\Om_\theta}$ can be expressed for $x \in \partial B$ as 
$$
    n_{\Om_\theta}\left(\phi_\theta(x)\right) = \frac{\left(D\phi_\theta\right)^{-T}(x) n_B(x)}{\left|\left(D\phi_\theta\right)^{-T}(x) n_B(x)\right|},
$$
and the change of variable formula leads to 
$$
    \int_{\dOm_\theta} H_{\Om_\theta}^2(y) dy = \int_{\partial B} H^2_{\Om_\theta}(\phi_\theta(x)) \Jac_{\partial B}(\phi_\theta(x)) dx.
$$

\subsection{PDE-related quantities}

In the sequel, we will consider several PDE related quantities, namely the torsion and Laplace eigenvalues of a domain. The torsion of a domain $\Om$ is defined as 
$$
    T(\Om) = \int_\Om udx ,
$$
where $u \in H^1_0(\Om)$ is the solution to
\begin{equation}
    \begin{cases}
        -\Delta u &=\ \  1\ \ \ \mbox{ in } \Om,\\
        \ \ \ \ \ u &=\ \  0\ \ \   \mbox{ on } \dOm.
    \end{cases}
\end{equation}
Equivalently, by taking $\phi(x) = u(x) + \frac{x_1^2}{2}$, we can look for a function $\phi$ which is solution to
\begin{equation}
    \begin{cases}
        -\Delta \phi &= \ \ 0\ \ \ \  \mbox{ in } \Om,\\
        \ \ \ \ \ \phi &= \ \ \frac{x_1^2}{2}\ \ \    \mbox{ on } \dOm.
    \end{cases}
\end{equation}

One may then seek a solution $\phi$ expressed as a linear combination of fundamental solutions, namely,
$$
    \phi(x) := \sum_{i=1}^n c_i \psi(x - y_i),
$$
where $\psi$ is the fundamental solution to $-\Delta \psi = \delta_0$ in $\R^d - \{0\}$ and $y_1,\dots, y_n \in \Om^c$. Since $\phi$ is a harmonic function in $\Om$, we only have to fit the boundary condition, for instance, in an $L^2$ sense, which amounts to solving
$$
    \min_{c_1, \dots, c_n \in \R} \int_{\dOm} \left| \phi(x) - \frac{x_1^2}{2} \right|^2d\sigma.
$$
Since this integral cannot be analytically computed for a general $\Om$, it is discretized and the resulting least squares problem is solved. This method takes inspiration from \cite{Hoskins2019Dec}, and in the Method of Fundamental Solutions that have already been used to solve shape optimization problems \cite{Antunes2022May, Bogosel2016Nov}.

A key advantage of the present method is that it allows to autodifferentiate the resulting $\theta\longmapsto T(\Omega_\theta)$ with respect to the parameters of the network, since every operation is supported by PyTorch (in particular, the solution of the least square problem via \verb|torch.linalg.lstsq|).

As previously mentioned, we will also consider the computation of Laplace eigenvalues with a homogeneous Neumann boundary condition i.e., positive values for which the equation
\begin{equation*}
    \begin{cases}
        -\Delta u &=\ \ \  \mu u\ \   \mbox{ in } \Om,\\
        \ \  \partial_nu &=\ \ \ \ \  0\ \     \mbox{ on } \dOm,
    \end{cases}
\end{equation*}
admits non-trivial solutions. 

It is well known that for a Lipschitz domain $\Om$, the spectral theorem guarantees that there exists an increasing sequence of eigenvalues 
$$
    0 = \mu_0(\Om) \leq \mu_1(\Om) \leq \dots \leq \mu_k(\Om) \leq \dots \to \infty.
$$
In order to numerically approximate the eigenvalues, we will use a mesh free RBF-Galerkin method as described in \cite{Wendland}. Let us quickly explain how this works in our case. Let
$$
    \int_{\Om_\theta} \nabla u \cdot  \nabla v dx = \mu \int_{\Om_\theta} uvdx,
$$
$u,v \in H^1(\Om)$ be the weak formulation of \cref{eq:neumann}. Using change of variables, we can pull back the equation over the reference domain $B$, which reads 
$$
    \int_{B} A_\theta \nabla u \cdot \nabla vdx = \mu \int_B (\Jac \phi_\theta) u vdx,
$$
for all $u,v \in H^1(B)$, where $A_\theta := (\Jac \phi_\theta) (D\phi_\theta)^{-1} (D\phi_\theta)^{-T}$. This problem is then discretized on a subspace spanned by radial basis functions (RBFs) $\varphi_i(x) := \psi(|x - x_i|)$, $x_i \in \Om$, leading to a system 
$$
    K_\theta \bar u = \mu M_\theta \bar u,
$$
where 
$$
    \left[K_\theta\right]_{ij} = \int_{B} A_\theta \nabla \varphi_i \cdot \nabla \varphi_j dx 
    \qquad \text{and} \qquad 
    \left[M_\theta\right]_{ij} = \int_B (\Jac \phi_\theta) \varphi_i \varphi_jdx.
$$

Since PyTorch does not currently offer automatic differentiation for generalized eigenvalue problems, we instead solve the equivalent problem
\begin{equation}
    \label{eq:ev}
     L_\theta^{-1} K_\theta L_\theta^{-T} \bar u = \mu \bar u,
\end{equation}
where $L_\theta$ is obtained via the Cholesky decomposition $M_\theta = L_\theta L_\theta^T$.

This method is chosen because it provides theoretical convergence guarantees (unlike other meshless methods relying on the strong formulation of the PDE like, for instance, the Kansa method \cite{fasshauer2007meshfree}, or PINN-based methods) and is easy to implement in PyTorch. An alternative, robust approach would consist in using finite element methods. Such a strategy would, however, require additional care to ensure a consistent interaction between the finite element solver and PyTorch, particularly with regard to the integration of automatic differentiation.

\subsection{A note on the practical computation of shape quantities}

In practice, all Jacobians can be computed using automatic differentiation. In this regard, our method is highly flexible, and can be readily adapted to a wide range of problems.

To compute integral quantities, one needs to choose some discretization of $B$ and $\partial B$. One might choose to use some Monte-Carlo approach, which has the advantage of being insensible to the curse of dimensionality. However, we will only consider convex sets in $\R^2$ and $\R^3$, and we will use an l-BFGS algorithm with line-search, which is very sensitive to noise. Hence, we decided to discretize $B$ by lattice points, while the sphere $\partial B$ is sampled via a Fibonacci lattice \cite{Gonzalez2010Jan}.

If the convex sets become too flat during the optimization, the computation of geometric and spectral quantities often loses accuracy. In order to avoid this, we add a regularizer based on the Jacobian of $\phi_\theta$. More precisely, we add a regularizer based on the maximal condition number of the
Jacobian over $B$. The condition number of a matrix $A \in \R^{d\times d}$, denoted $\text{cond} A$, is the ratio between its largest and smallest singular value. The regularizer then reads
$$
    \sum_{i=1}^N \| \text{cond} \Jac \phi_{\theta_i} \|_\infty.
$$
Hence, the total loss function that we minimize is
$$
    L(\theta_1, \dots, \theta_N) = \sum_{1 \leq i\neq j \leq N} \left|F(\Om_{\theta_i}) - F(\Om_{\theta_j})\right|^{-s} + \alpha \sum_{i=1}^N \| \text{cond} \Jac \phi_{\theta_i} \|_\infty.
$$

\subsection{Code availability}

All these functions have been fully implemented in PyTorch and integrated into an easy-to-use package, available at \url{https://github.com/EloiMartinet/ConvexShapeOpt}. Detailed installation instructions and comprehensive documentation are provided to assist the reader in conducting their own numerical experiments.

To facilitate reproducibility and ease of use, a collection of notebooks is available \href{https://drive.google.com/drive/folders/12-tWL7n0nP48uCksgxboSXp1QX3xLcwN}{here}, allowing the reader to experiment with the code without any prior installation. The repository contains three scripts that install the package and solve representative shape optimization problems.

\section{Obtained results}

\subsection{The diagram \texorpdfstring{$\Vol$, $\Per$, $\W$}{Vol, Per, W}}

In this section, we are interested in the diagram of the volume, the perimeter and the moment of inertia, see Section \ref{ss:integral_quantities} for the definitions. Our method is applied to obtain numerical descriptions for three classes of sets: the planar convex domains, the planar convex domains with two orthogonal axes of symmetries a class that has been studied in detail in \cite{Gastaldello2023Jul}, and the convex domains in three dimensions.  

A first step is to reduce the three quantities into two scale invariant ones. Using a change of variables, we can easily show that
$$\Vol(t\Om) = t^d \Vol(\Om), \quad \Per(t\Om) = t^{d-1} \Per(\Om), \quad \text{and}\quad \W(t\Om) = t^{d+2} \W(\Om),$$
which means that the quantities
$$\frac{\Vol(\Om)^{(d+2)/d}}{\W(\Om)} \quad \mbox{ and } \quad \frac{\Vol(\Om)}{\Per(\Om)^{d/(d-1)}}$$
are scale-invariant. Moreover, the isoperimetric inequality implies that
$$\frac{\Vol(\Om)}{\Per(\Om)^{d/(d-1)}} \leq \frac{\Vol(B)}{\Per(B)^{d/(d-1)}}.$$

One can apply the Hardy-Littlewood inequality for symmetric rearrangements to show that
$$
    \frac{\Vol(\Om)^{(d+2)/d}}{\W(\Om)} \leq \frac{\Vol(B)^{(d+2)/d}}{\W(B)},
$$
where $B$ is any ball (see \cite[Chapter 7, Section 2]{Polya1951}). Hence, by taking
$$
    x(\Om) = \frac{\Vol(\Om)^{(d+2)/d}}{\W(\Om)} \frac{\W(B)}{\Vol(B)^{(d+2)/d}}
    \quad \mbox{ and } \quad
    y(\Om) = \frac{\Vol(\Om)}{\Per(\Om)^{d/(d-1)}} \frac{\Per(B)^{d/(d-1)}}{\Vol(B)},
$$
we have that the diagram
$$
    \mathcal{D}^d_{VPW} := \left\{ \left(x(\Om), y(\Om)\right) : \Om \in \mathcal{K}_d\right\}
$$
is a subset of the unit square $[0,1]^2$. 

\begin{remark}
    Rescaling the problems so that the diagram fits within $[0,1]^2$ requires additional information about the shape functions involved, but this is by no means required by the method itself. Nevertheless, for the sake of consistency, we will adopt this rescaling throughout what follows.
\end{remark}

\subsubsection{The diagram \texorpdfstring{$\Vol$, $\Per$, $\W$}{Vol, Per, W} for planar convex bodies}

When considering planar convex sets $d=2$, we take 
$$
    x(\Om) = \frac{1}{2\pi}\frac{\Vol(\Om)^2}{\W(\Om)}
    \quad \mbox{ and } \quad
    y(\Om) = 4\pi \frac{\Vol(\Om)}{\Per(\Om)^2},
$$
In this case, we recall the classic inequality obtained by G. Pólya in \cite[Theorem 1]{zbMATH03175801} 
$$\frac{\Vol(\Omega)\Per(\Omega)^2}{W(\Omega)}> 48,$$
which is sharp as the equality asymptotically holds when considering a sequence of thinning stadiums. This inequality is then written in terms of $x(\Omega)$ and $y(\Omega)$ as follows 
$$y(\Omega)<\frac{\pi^2}{6}x(\Omega).$$

Also, we recall the classic conjecture by G. Pólya that states that for all planar convex sets
$$\frac{\Vol(\Omega)\Per(\Omega)^2}{W(\Omega)}\leq 108,$$
with equality if and only if $\Omega$ is an equilateral triangle. This estimate is written in terms of $x(\Omega)$ and $y(\Omega)$ as follows:
$$y(\Omega)\ge\frac{2\pi^2}{27}x(\Omega).$$

Note that this diagram has been studied recently in \cite{Gastaldello2023Jul} in the planar case $d=2$ for convex bodies with two orthogonal axes of symmetry, which is addressed in Section \ref{sss:doubly_sym}.

We now solve \cref{eq:main-opt} for $F(\Om) := \left(x(\Om), y(\Om)\right)$ where $\Om \subset \mathcal{K}_d$, with $d\in\{2,3\}$. Note that while the $\R^2$ case can still run on a laptop, the $\R^3$ case runs in a few hours on an Nvidia h100 GPU. This is due to the fact that the number of quadrature points necessary to compute the quantities of interest to the desired accuracy in $\R^3$ is significantly larger than in $\R^2$, preventing the computational graph of PyTorch to fit in a standard laptop RAM.

\begin{figure}[h!]
    \centering
    \includegraphics[width=0.8\linewidth]{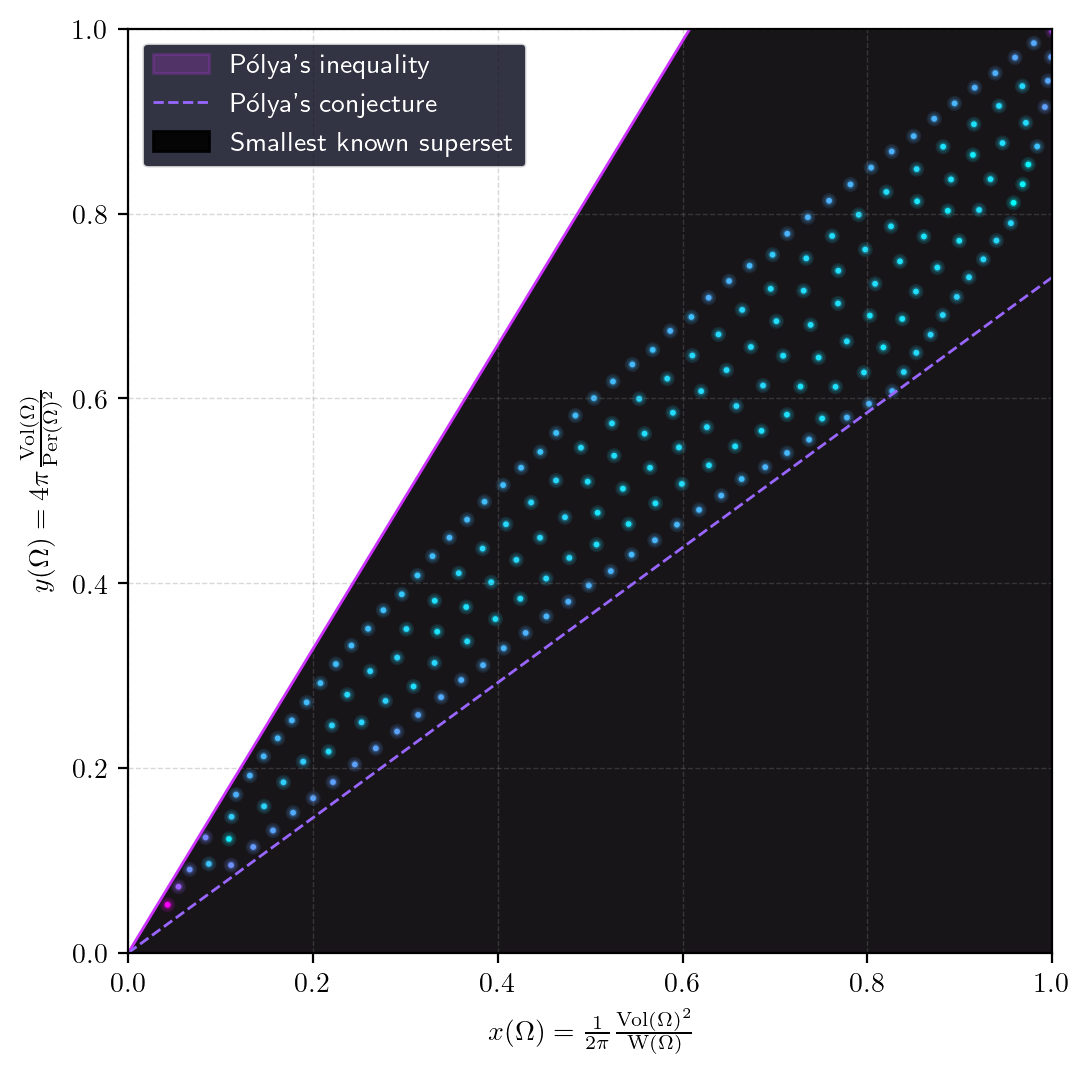}
    \caption{Numerical approximation of the diagram  $\mathcal{D}^2_{VPW}$ for $d=2$.}
    \label{fig:APW_2d}
\end{figure}

The resulting diagram is displayed in \cref{fig:APW_2d}. The black region represents the smallest known set that contains the diagram. The blue points are the optimal points given by our method. As you can see, the points are uniformly distributed in the diagram.

A visualization of the sets comprising the diagram can be found \href{https://eloimartinet.github.io/projects/blaschke-santalo/diagrams/plot_diagram.html?folder=APW_2d/}{here}. As observed, the neural network representation of the convex shapes allows for both smooth and polygonal shapes. The upper part of the boundary seems to be realized by smooth domains resembling stadiums, while the leftmost lower part is apparently realized by isosceles triangles. We see that no points are sampled near the origin. The reason is the Jacobian regularizer preventing flattening sets. Indeed, a small $y$ indicates a large isoperimetric ratio, which can only be attained for very flat convex sets.

\subsubsection{The diagram \texorpdfstring{$\Vol$, $\Per$, $\W$}{Vol, Per, W} for doubly symmetric planar convex bodies}\label{sss:doubly_sym}

\begin{figure}
    \centering
    \includegraphics[width=0.8\linewidth]{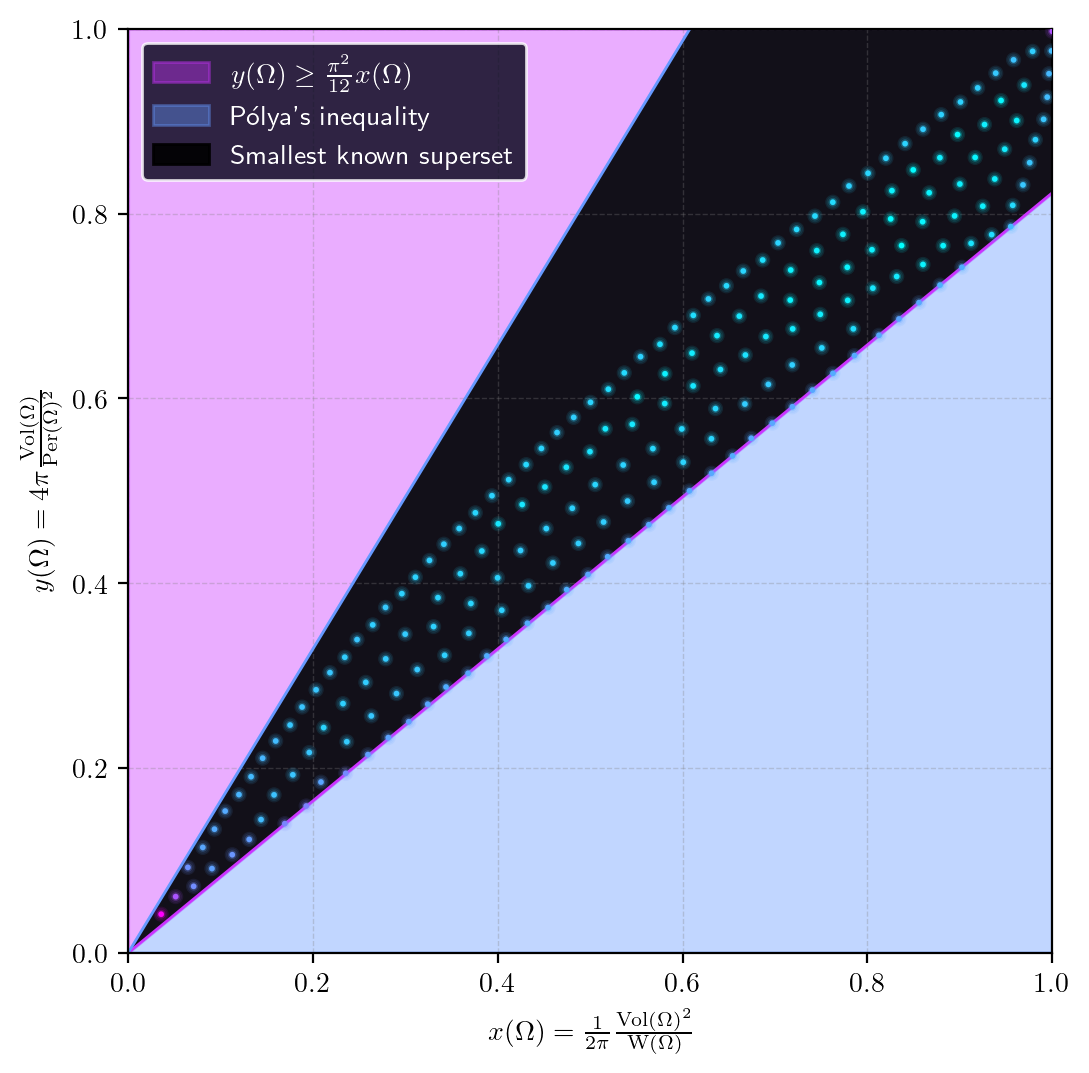}
    \caption{Numerical approximation of the diagram  $\mathcal{D}^2_{VPW}$ for planar and doubly symmetric convex bodies, along with the theoretically known lower boundary.}
    \label{fig:APW_2d_sym}
\end{figure}

We check that our method correctly samples the diagram studied in \cite{Gastaldello2023Jul}, where planar convex bodies with two orthogonal axes of symmetry are considered. In the cited paper, a large part of the lower boundary is explicitly described (see \cite[Corollary 5.2]{Gastaldello2023Jul}) by the inequality
$$
    y(\Om) \geq \frac{\pi^2}{12} x(\Om).
$$
By taking the subgroup of isometries of $\R^2$
$$
    G := \left\{
        (x,y) \mapsto (x,y), \quad
        (x,y) \mapsto (-x,y),\quad
        (x,y) \mapsto (x,-y),\quad
        (x,y) \mapsto (-x,-y)
    \right\}
$$
in \cref{eq:symmetric}, we are able to represent only convex sets with horizontal and vertical reflection symmetries. The resulting diagram is displayed in \cref{fig:APW_2d_sym}. 


A visualization of the sets comprising the diagram can be found \href{https://eloimartinet.github.io/projects/blaschke-santalo/diagrams/plot_diagram.html?folder=APW_2d_sym/}{here}. We observe that our method accurately recovers the known part of lower boundary of the diagram, which is attained by rhombi, see \cite[Theorem 5.1]{Gastaldello2023Jul}. On the other hand, the (symmetric) sets located on the upper boundary seem to be the same as the ones obtained for the general case presented in fig. \ref{fig:APW_2d}.

\subsubsection{The diagram \texorpdfstring{$\Vol$, $\Per$, $\W$}{Vol, Per, W} for convex bodies in three dimensions}

\begin{figure}[h!]
    \centering
    \includegraphics[width=0.8\linewidth]{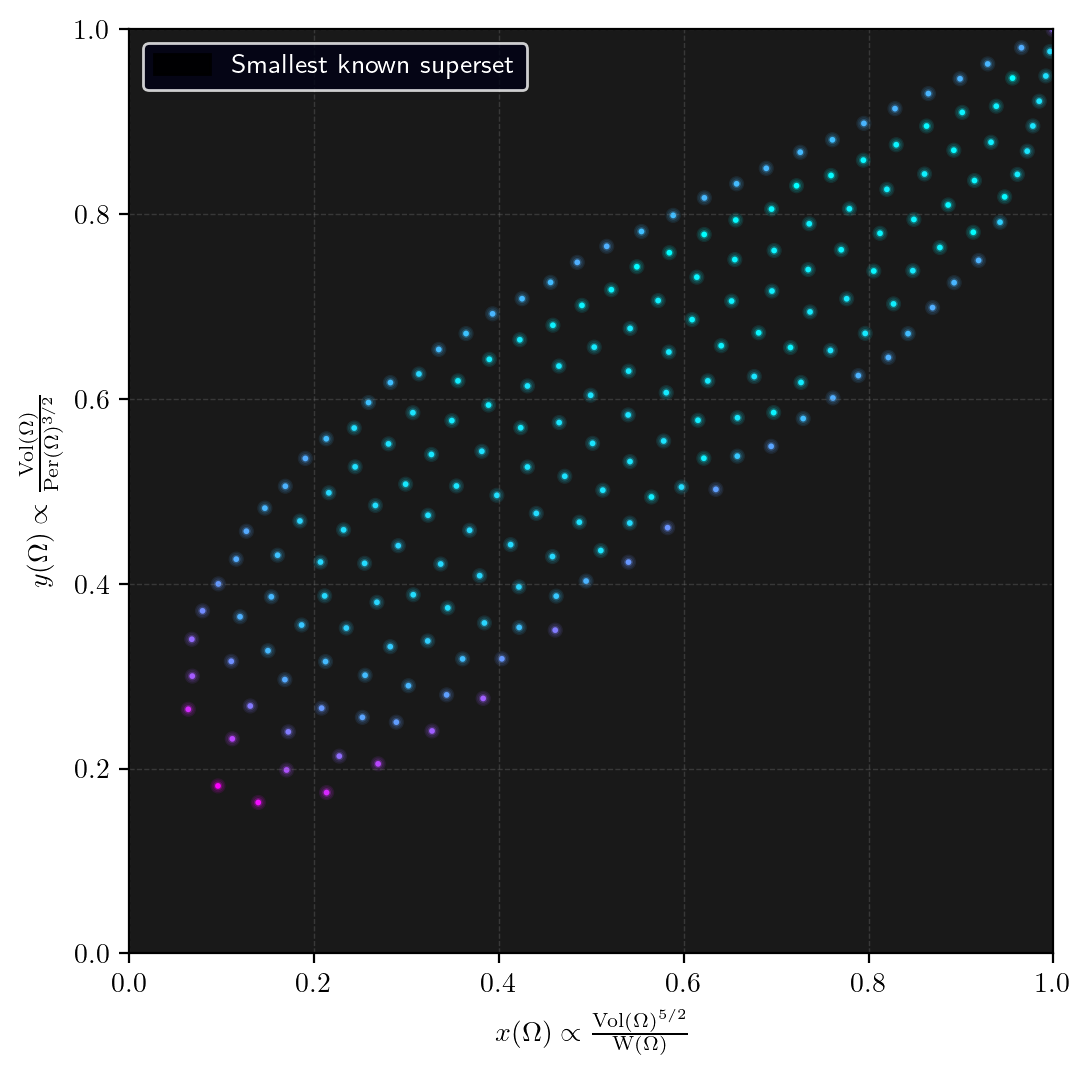}
    \caption{Numerical approximation of the diagram $\mathcal{D}^3_{VPW}$ for $d=3$}
    \label{fig:APW_3d}
\end{figure}

We show that the method also allow for diagrams involving 3 dimensional convex sets. In this case,
$$
    x(\Om) = \left(\frac{3}{4}\right)^{5/3} \left(\frac{4}{5\pi^{2/3}}\right) \frac{\Vol(\Om)^{5/3}}{\W(\Om)} 
    \quad \mbox{ and } \quad
    y(\Om) = 6 \sqrt{\pi}\frac{\Vol(\Om)}{\Per(\Om)^{3/2}}.
$$
The resulting diagram is displayed in \cref{fig:APW_3d}. A visualization of the sets comprising the diagram can be found \href{https://eloimartinet.github.io/projects/blaschke-santalo/diagrams/plot_diagram.html?folder=APW_3d/}{here}. The upper part of the boundary appears to be attained by sausage-shaped domains, namely domains obtained as revolutions of stadiums about their axes of symmetry. The leftmost lower portion is apparently attained by flattening polyhedra.

\subsection{The diagram $\Vol$, $\Per$, $\mathbf{T}$}

To illustrate the applicability of our method to more complex PDE-dependent problems, we consider a diagram involving the torsional rigidity functional. Let us recall that the torsional rigidity of a domain $\Om$ is given by
$$
    T(\Om) = \int_\Om u dx,
$$
where $u \in H^1_0(\Om)$ is the solution of
\begin{equation}
    \begin{cases}
        -\Delta u &=\ \ 1 \ \mbox{ in } \Om,\\
        \ \ \ \ \ u &=\ \ 0 \   \mbox{ on } \dOm.
    \end{cases}
\end{equation}

A straightforward change of variables shows that $T(t\Om) = t^{d+2}T(\Om)$. In order to obtain scale-invariant quantities and to have a diagram contained in the unit square $[0,1]^2$, we set
$$
    x(\Om) := \frac{\Vol(\Om)^{(d-1)/d}}{\Per(\Om)}\cdot \frac{\Per(B)}{\Vol(B)^{(d-1)/2}}  \quad \quad \text{and} \quad \quad y(\Om)= \frac{\Vol(B)^{(d+2)/d}}{T(B)}\cdot  \frac{T(\Om)}{\Vol(\Om)^{(d+2)/d}},
$$
then, consider the diagram $\mathcal{D}_{VPT}:=\{(x(\Omega),y(\Omega)):\ \Omega\in \mathcal{K}_d\}$. 
The fact that $y(\Om) \in [0,1]$ is a consequence of the classic Saint-Venant inequality \cite{de1856memoire}.

Note that, by separation of variables, the torsion of the unit ball $B$ can be computed explicitly and is then given by
$T(B) = \omega_d/(d(d+2))$, where $\omega_d$ is the Lebesgue measure of $B$.

When $d=2$, Makai \cite{Makai} also proved the following inequality:
$$
    \frac{T(\Om)\Per(\Om)^2}{\Vol(\Om)} \leq \frac{2}{3},
$$
while Pólya \cite{Polya} showed that
$$
    \frac{T(\Om)\Per(\Om)^2}{\Vol(\Om)} \geq \frac{1}{3}.
$$

In terms of $x$ and $y$, those inequalities translate to
$$
    \frac{2}{3} x(\Om)^2 \leq y(\Om) \leq \frac{4}{3} x(\Om)^2,
$$
which yields the inclusion 
$$\mathcal{D}^2_{VPT}\subset  \left\{(x,y):\ \frac{2}{3} x^2 \leq y \leq \frac{4}{3} x^2\right\}.$$

\begin{figure}[h!]
    \centering
    \includegraphics[width=0.8\linewidth]{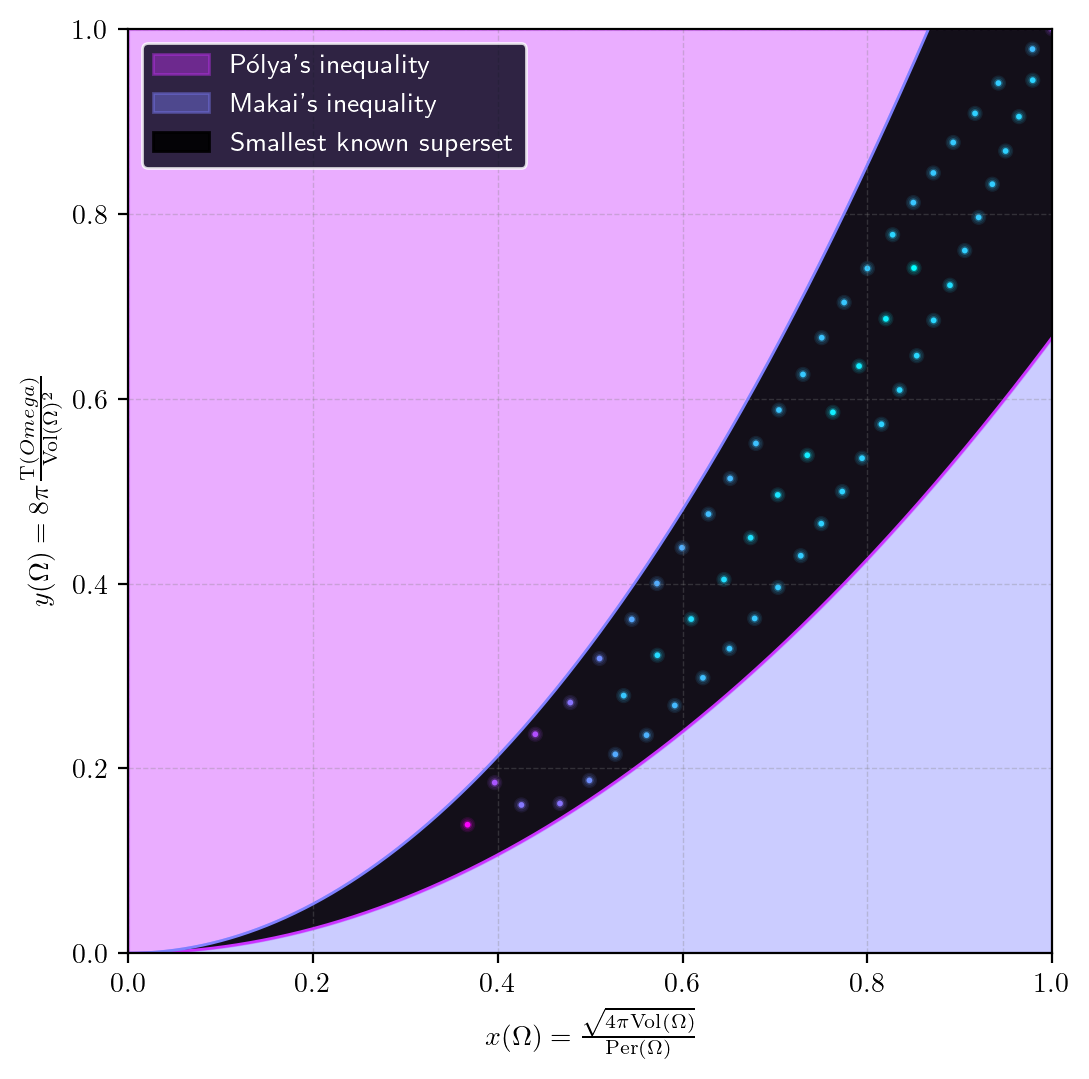}
    \caption{Numerical approximation of the diagram  $\mathcal{D}^2_{VPT}$ for $d=2$ along with the smaller theoretically known  superset.}
    \label{fig:APT}
\end{figure}

The resulting diagram is displayed in \cref{fig:APT}. A visualization of the sets comprising the diagram can be found \href{https://eloimartinet.github.io/projects/blaschke-santalo/diagrams/plot_diagram.html?folder=APT_2d/}{here}.




\subsection{The diagram \texorpdfstring{$\Vol$, $\mathbf{\mu_1}$, $\mathbf{\mu_2}$}{Vol, mu1, mu2}}

In this section, we study the diagram of the volume and the first two non-trivial Neumann eigenvalues, which are positive $\mu$ for which the problem
\begin{equation}
    \label{eq:neumann}
    \begin{cases}
        -\Delta u &=\ \ \  \mu u\ \   \mbox{ in } \Om,\\
        \ \  \partial_nu &=\ \ \ \ \  0\ \     \mbox{ on } \dOm,
    \end{cases}
\end{equation}
admits non-trivial solutions. We note that this diagram has been studied by P. Antunes and A. Henrot in \cite{Antunes2010}. 

It can be shown that for every $k \in \N$, the function $t \mapsto \mu_k(t\Om)$ is $(-2)$-homogeneous, which implies that the quantity 
$$
    \Vol(\Om)^{2/d}\mu_k(\Om)
$$
is scale-invariant for every $k \in \N$. Using Brouwer's fixed point argument and properties of the eigenfunctions of the ball, Weinberger showed that for $k=1$ and every $d \in \N^*$, the previous quantity is maximized by the ball \cite{weinberger1956isoperimetric}. Using a more involved argument of topological degree theory, D. Bucur and A. Henrot \cite{Bucur2019Jun} demonstrated that the previous quantity is maximized by the union of two disjoint balls of the same measure, for $k=2$ and all $d \in \N^*$.

For $d=2$, we consider the quantities 
$$
    x(\Om) := \frac{\Vol(\Om)\mu_1(\Om)}{\Vol(B)\mu_1(B)} \qquad \text{and} \qquad  y(\Om) := \frac{\Vol(\Om)\mu_2(\Om)}{\Vol(B \sqcup B)\mu_2(B \sqcup B)},
$$
where $B \sqcup B$ here denotes the union of two disjoint balls of the same measure.

The Weinberger and Bucur--Henrot inequalities translate, respectively, into $x(\Om) \leq 1$ and $y(\Om) \leq 1$, which means that $$\mathcal{D}^2_{V\mu} := \{(x(\Omega),y(\Omega))\ :\ \Omega\in \mathcal{K}_2\} \subset [0,1]^2.$$ 

Note that since we only consider convex sets, we have $y(\Om) < 1$, the optimal value being unknown. We refer to \cite[Section 3]{Antunes2010} for numerical results on the maximizer of $\Vol(\Omega)\mu_2(\Omega)$  among planar convex sets. Because of the ordering of the eigenvalues, we also have $\Vol(\Om) \mu_1(\Om) \leq \Vol(\Om) \mu_2(\Om)$. Since $\Vol(B \sqcup B)\mu_2(B \sqcup B) = 2\Vol(B)\mu_1(B)$, this leads to the inequality
$$
    x(\Om) \leq 2 y(\Om).
$$

Meanwhile, Szegö \cite{zbMATH03087123} showed that for planar and simply connected domains, 
$$
    \frac{1}{\Vol(\Om)\mu_1(\Om)} + \frac{1}{\Vol(\Om)\mu_2(\Om)} \geq \frac{1}{\Vol(B)\mu_1(B)} + \frac{1}{\Vol(B)\mu_2(B)},
$$
which translates, in terms of $x$ and $y$, to 
$$
    y(\Om) \leq \frac{x(\Om)}{4 x(\Om) - 2}.
$$

\begin{remark}
    In a recent paper \cite{Funano2025Nov}, K. Funano proved that for every convex body in $\mathbb{R}^d$ the following estimates relating $\mu_k$ and $\mu_{k+1}$
$$\mu_{k+1}(\Omega)\leq C d^7 \mu_{k}(\Omega),$$
where $C$ is a universal constant independent of $k$ and $d$. In the same vein, he proved for planar convex bodies the following result
$$\mu_k(\Omega)\leq C' \left(\frac{k}{\ell} \right)^2 \mu_\ell(\Omega),$$
where $k\ge \ell$ and $C'$ is a universal constant. When taking $k=2$ and $\ell=1$, this result supports the conjecture $\mu_2(\Omega)\leq 4 \mu_1(\Omega)$ stated in \cite{zbMATH00238628} and numerically observed in \cite{Antunes2010} and in the present work. 
\end{remark}





\begin{figure}[h!]
    \centering
    \includegraphics[width=0.8\linewidth]{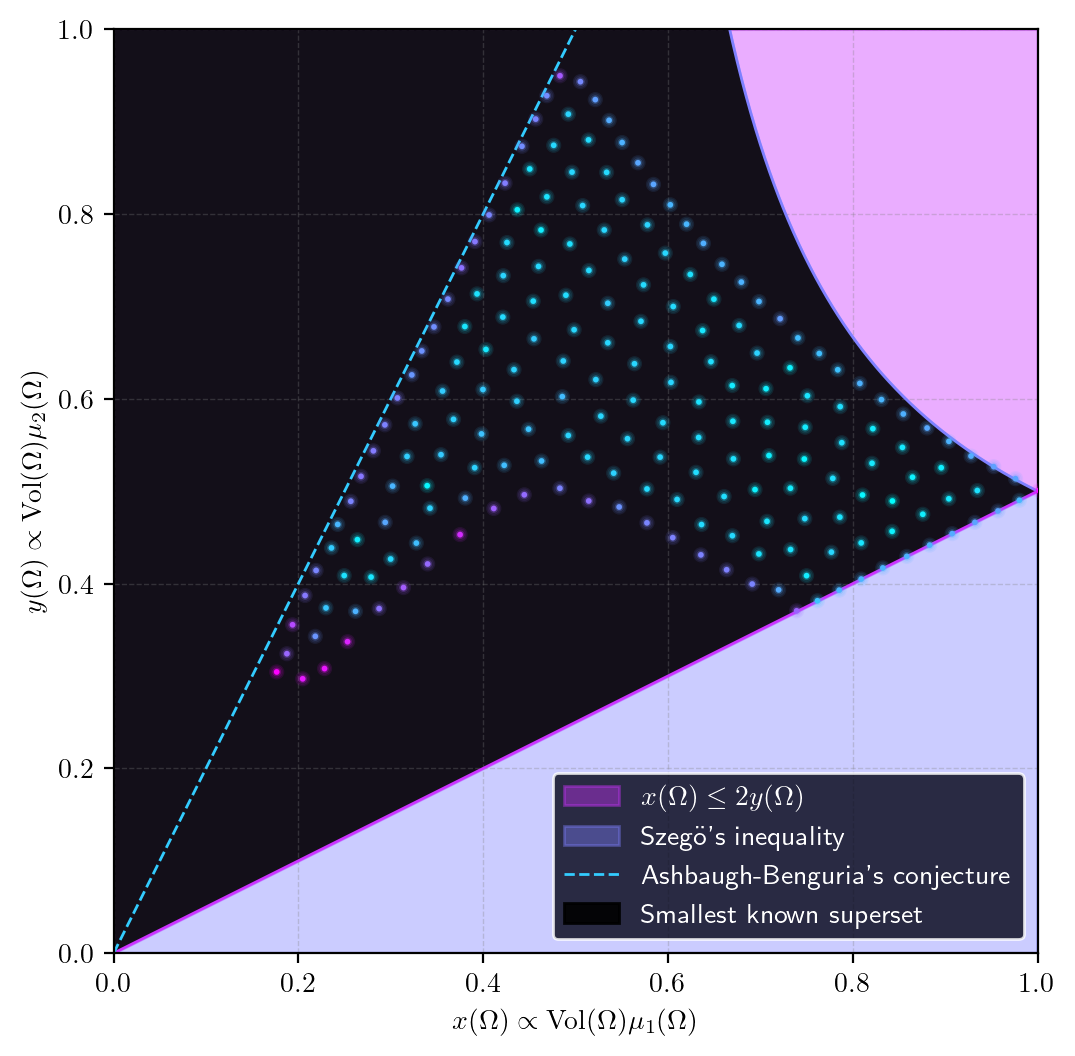}
    \caption{Numerical approximation of the diagram $\mathcal{D}^2_{V\mu}$ for $d=2$ along with the theoretically known smaller superset.}
    \label{fig:Neumann_2d}
\end{figure}

In fig. \ref{fig:Neumann_2d}, we present the obtained diagram for the class of planar convex sets. The shapes can be visualized \href{https://eloimartinet.github.io/projects/blaschke-santalo/diagrams/plot_diagram.html?folder=Neumann_2d}{here}.

 In contrast to the other diagrams, we are not able to recover with high precision the extremal shapes corresponding the the upper part of the boundary of the diagram. In particular, the shapes of the top left part of the boundary appears to be rounded rectangles, although actual rectangles yield a better value (which agrees with the Ashbaugh-Benguria's conjecture). This is due to the fact that because of the approach that we chose, PyTorch needs to differentiate through the whole assembly and the computation of the eigenvalues of ~\cref{eq:ev}, which requires a large amount of memory and limits the quadrature precision. A natural direction for future research is to enhance the memory efficiency of the proposed approach. 
 
\subsection{The diagram \texorpdfstring{$\Vol$, $\Per$, $\mathbf{E}$}{Vol, Per, E}}

In this section, we show that our approach can be easily adapted to problems involving higher-order boundary terms, such as the curvature of the boundary. Here, we consider the Willmore energy (also called \textit{elasticae} for the 1 dimensional boundary of planar sets) defined as 
$$
    E(\Om) := \int_\dOm H_\Om^2d\sigma,
$$
where $\Omega$ is a $C^2$ convex body of $\mathbb{R}^d$ and $H_\Om$ is the mean curvature of $\dOm$. We note that the functional $E$ is $(d-3)$-homogeneous and we denote by $\widetilde{\mathcal{K}}_d$ the class of $C^2$ convex bodies of $\mathbb{R}^d$.

\subsubsection{The diagram \texorpdfstring{$\Vol$, $\Per$, $E$ for planar convex sets}{Vol, Per, E}}

In the planar case $d=2$, the functional $E$ is $(-1)$-homogeneous. We then consider the following scale invariant quantities
$$
    x(\Om) := \frac{4 \pi \Vol(\Om)}{\Per(\Om)^2}
    \qquad \text{and} \qquad 
    y(\Om) := \frac{2 \pi^2}{\Per(\Om) E(\Om)},
$$
and the diagram $\D^2_{VPE}:=\{(x(\Omega),y(\Omega)):\ \Omega\in \widetilde{\mathcal{K}}_2\}$, which was studied in \cite{BHT}, where qualitative results concerning the diagram and the extremal shapes are established. 

The isoperimetric inequality ensures that $x(\Om) \leq 1$. On the other hand, since
$$
    2 \pi = \int_\dOm H_\Om \leq \left( \int_\dOm H_\Om^2 \right)^\frac{1}{2} \left( \int_\dOm 1 \right)^\frac{1}{2} = E(\Om)^\frac{1}{2} \Per(\Om)^\frac{1}{2},
$$
we have $y(\Om) \leq 1$. At last, an important inequality, due to Gage \cite{Gage1983Dec}, states that 
$$
    \frac{E(\Om) \Vol(\Om)}{\Per(\Om)} \geq \frac{\pi}{2},
$$
which translates, in our case, to the inequality
$y(\Om) \leq x(\Om)$.


\begin{figure}[!h]
    \centering
    \includegraphics[width=0.8\linewidth]{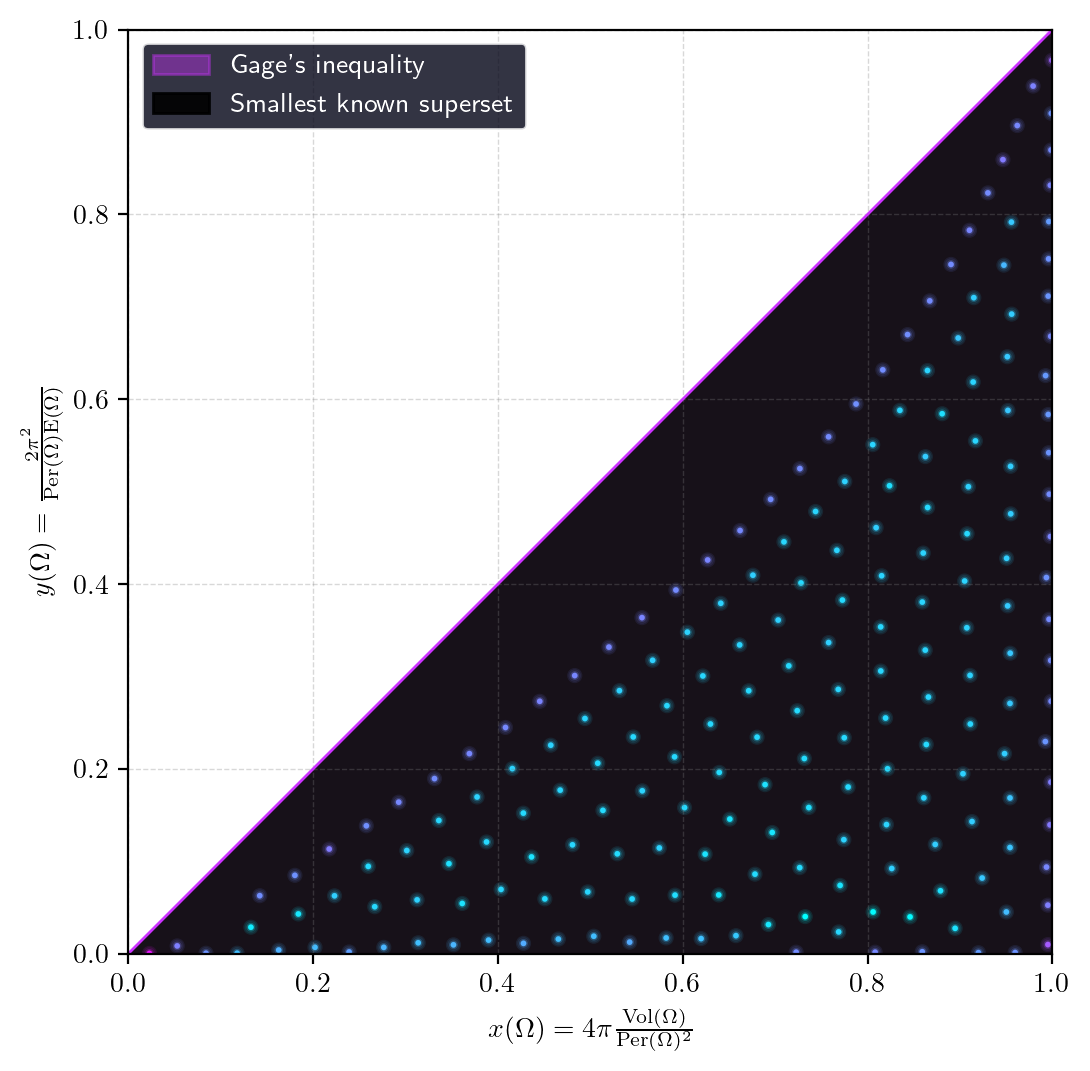}
    \caption{Numerical approximation of the diagram  $\mathcal{D}^2_{VPE}$ for $d=2$ along with the theoretically known smaller superset.}
    \label{fig:VPE_2d}
\end{figure}

The resulting diagram is displayed in \cref{fig:VPE_2d}. A visualization of the sets comprising the diagram can be found \href{https://eloimartinet.github.io/projects/blaschke-santalo/diagrams/plot_diagram.html?folder=VPE_2d/}{here}.

The lower part of the boundary of the diagram $\{(x,0): x\in [0,1]\}$ is approached by rounded polygons i.e., sets of the form $P_\eps =  P\oplus B_\eps(x)$, where $\oplus$ stands for the Minkowski sum $P$ is a convex polygon and $B_\eps(x)$ is a ball of center $x$ and radius $\eps$. Indeed, for any $x^* \in (0,1)$, there exists a polygon $P$ such that $x(P) = x^*$. Then, one can show that $x(P_\eps) \xrightarrow[\eps \to 0]{} x^*$ while $E(P_\eps)$ scales as $\eps^{-1}$. 

For the right part of the boundary i.e., $\{(0,y): y\in [0,1]\}$, one can approach any point $(0,y^*) \in [0,1]^2$ by taking sets obtained via Minkowski sums of regular polygons and the unit ball. More precisely, for $N\ge 3$, and $t\in(0,1]$, we take $R_N^t:=(1-t)R_N \oplus tB$, where $R_N$ is the regular polygon of $N$ sides and $B$ is a ball. Finally, it is straightforward that the continuous paths $\Gamma_N:=\{x(R_N^t),y(R_N^t): t\in(0,1]\}\subset \D^2_{VPE}$ converge to the segment  $\{(0,y): y\in [0,1]\}$ when $N$ tends to infinity.  

\subsubsection{The diagram \texorpdfstring{$\Vol$, $\Per$, $E$ for convex sets of dimension three}{Vol, Per, E}}

In dimension three, the Willmore energy is already a scale-invariant functional. Hence, we take 
$$
    x(\Om) = \frac{36 \pi \Vol(\Om)^2}{\Per(\Om)^3}
    \qquad \text{and} \qquad 
    y(\Om) = \frac{4 \pi}{E(\Om)},
$$
and consider the diagram $$\D^3_{VPE}:=\{(x(\Omega),y(\Omega)):\ \Omega\in \widetilde{\mathcal{K}}_3\}\subset [0,1]^2,$$ 
where, the normalization in $y$ is motivated by Willmore's inequality 
$$
    E(\Om) \geq 4 \pi,
$$
see for instance \cite{Agostiniani2020Feb} and the references therein.

\begin{figure}[!h]
    \centering
    \includegraphics[width=0.8\linewidth]{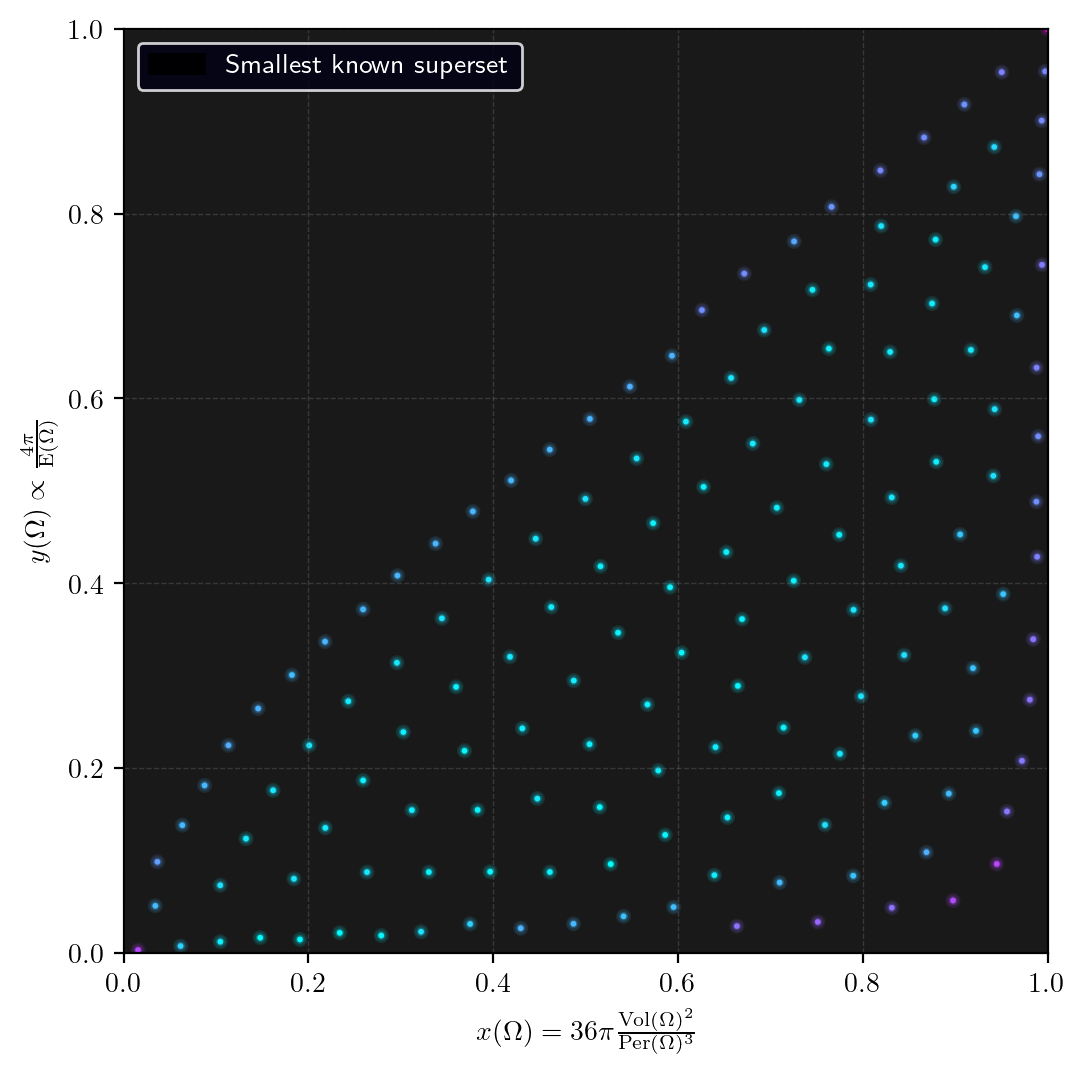}
    \caption{Numerical approximation of the diagram  $\mathcal{D}^3_{VPE}$ for $d=3$.}
    \label{fig:VPE_3d}
\end{figure}

The resulting diagram is displayed in \cref{fig:VPE_3d}. A visualization of the sets comprising the diagram can be found \href{https://eloimartinet.github.io/projects/blaschke-santalo/diagrams/plot_diagram.html?folder=VPE_3d/}{here}.

\subsection{Limitations of the method}

While we have seen through the previous examples that this method is both flexible and robust, it comes with some limitations. First, flexibility and ease of use comes with a price in terms of computation and memory efficiency; indeed, relying on automatic differentiation requires PyTorch to store a large computational graph, particularly when high precision necessitates the use of a substantial number of quadrature points. One potential idea for improvement would be to use higher-order quadrature schemes.

Another limitation comes from the representation itself of the shapes. Since we decided to use a representation \textit{via} gauge functions, some natural geometric quantities associated with convex sets (like the minimal width or the diameter) are harder to compute than when using support functions. However, we note that our proposed parametrization can also be used to parametrize support functions, as every sublinear function $p$ is the support function of the convex set $K = \{y \in \R^d : x \cdot y \leq p(x) \text{ for all } x \in \R^d \}$. However, the computation of quantities relying on change of variables (like the Neumann spectrum) might become more challenging.


\subsection{Code availability}

For reproducibility purposes, the code used to generate this figures is made fully accessible here: \url{https://github.com/EloiMartinet/ConvexShapeOpt}. The codes for the generations of the diagrams, listed in the \verb|scripts/| folder, can all run in less than $24$ hours on a Nvidia h100 GPU. However, most of them  (namley, the purely geometrical ones in two dimensions) needs much less and can also run on CPU in a few hours.

\section{Comparison with random sampling}

Generating random convex polygons (or, more generally, random convex sets) is a problem of intrinsic interest. In the present work, we employ an algorithm described \href{https://cglab.ca/~sander/misc/ConvexGeneration/convex.html}{here}, which builds upon the work of P. Valtr \cite{random_polygon}, where the author computes the probability that a set of $n$ points, independently and uniformly distributed in a given square, is in convex position. As noted in Section 4 of that work, the proof also leads to a fast and "uniform" procedure for generating random convex sets within a prescribed square.

\begin{figure}[!h]
    \centering
    \includegraphics[width=0.42\linewidth]{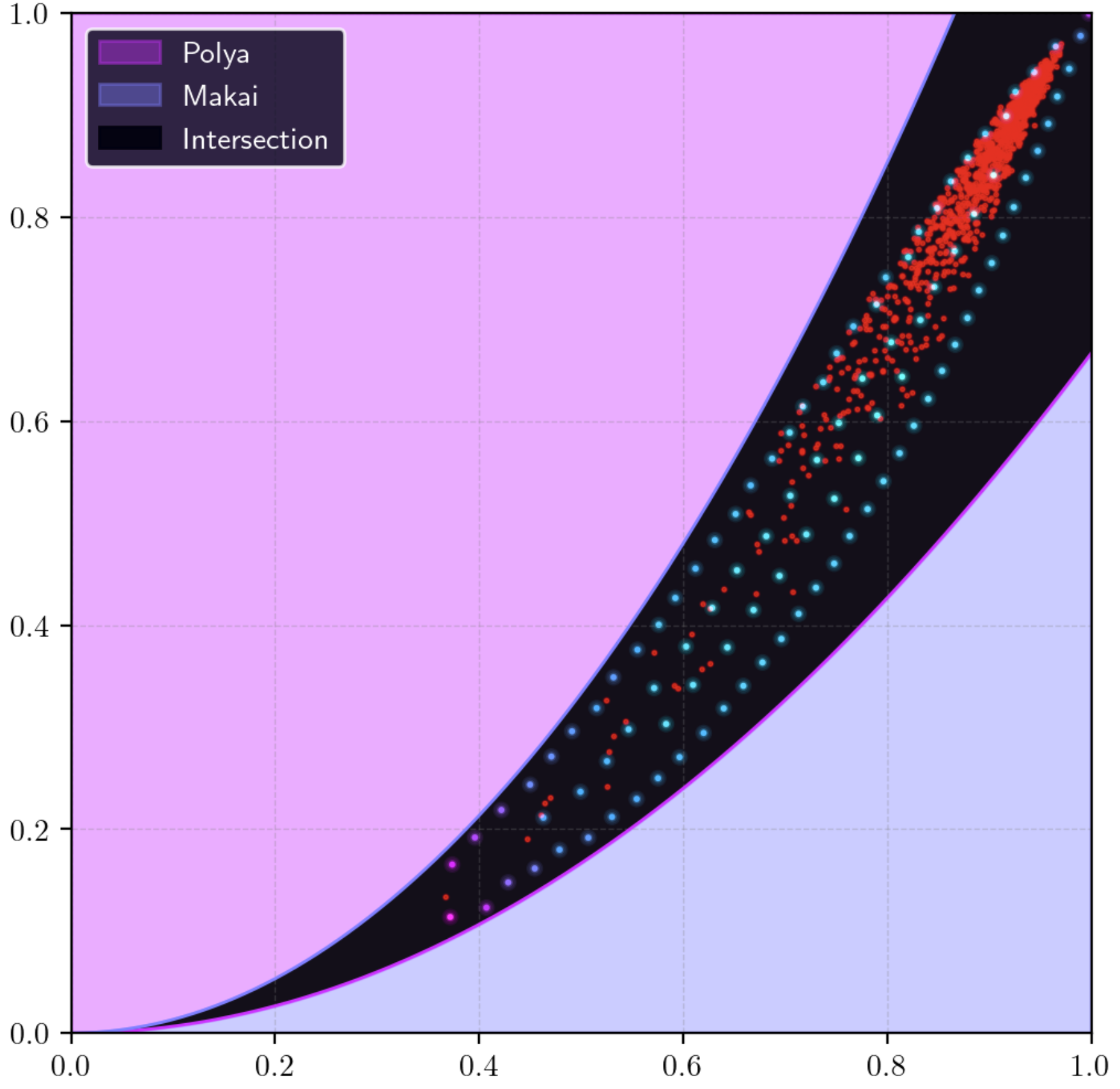}
    \includegraphics[width=0.44\linewidth]{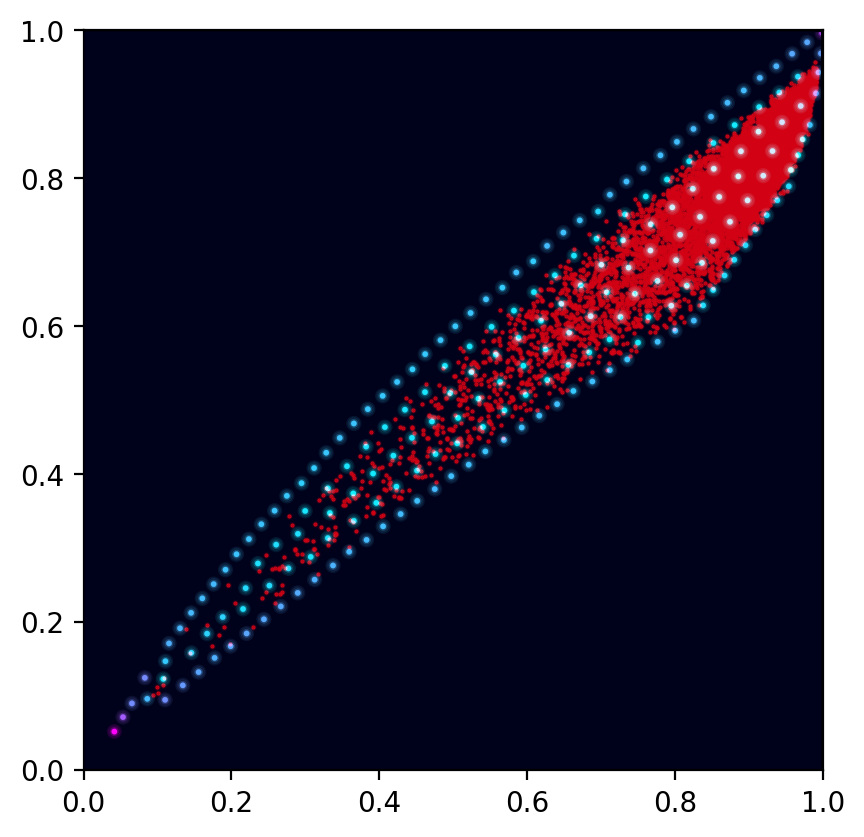}
    \caption{Comparison between our results and the ones obtained by random sampling of convex polygons (red dots) for the the diagram $\mathcal{D}^2_{VPW}$ on the left and the diagram $\mathcal{D}^2_{VPT}$ on the right.} 
    \label{fig:comparison}
\end{figure}

In fig. \ref{fig:comparison}, we compare our results with those obtained via the (naive) Monte Carlo–type method, which consists of generating random convex polygons, computing the associated functionals, and thereby producing a rudimentary approximation of the Blaschke--Santaló diagrams. As can be seen on the figure, the Monte-Carlo method fails to capture certain parts of the diagrams.

\section{Conclusion and perspectives}

In this work, we have introduced a novel numerical framework for the exploration of Blaschke--Santaló diagrams. By leveraging an original invertible neural network architecture to parametrize convex bodies, we developed a fully unconstrained method that avoids handling the convexity constraint explicitly, while it allows to approximate any convex set. We also take advantage of the automatic differentiation capability of PyTorch to handle the complex multi-shape optimization problem that we must solve to uniformly sample the diagrams.

The developed methodology proved highly effective across a diverse range of shape functionals: including purely geometric and PDE-based ones such as the torsional rigidity and the eigenvalues of the Neumann Laplacian.  Compared to traditional random sampling of polygons, the proposed method, based on interacting particle systems, yields a substantially more uniform coverage of the diagrams. For the dimensional scalability, it is worth noting that while the computational cost increases in three dimensions due to quadrature requirements, the method remains robust, providing satisfactory descriptions of the diagrams, without any additional implementation difficulties. Lastly, by leveraging automatic differentiation via PyTorch, the framework provides a versatile, easily adaptable system for new functionals, with the complete implementation accessible to the research community.

While the current results focus on the class of convex bodies, the ideas covered in this paper can be extended in several directions:
\begin{itemize}
    \item \textbf{Alternative representations}: One could investigate alternative formulations of sublinear neural networks, for instance by modeling them as a smoothed maximum of the gauge functions of disks rather than of half-spaces. Moreover, combining multiple representations could improve the sampling of Blaschke--Santal\'o diagrams.
    \item \textbf{General Invertible Architectures}: Extending the parametrization beyond the sublinear functions used here to more general diffeomorphisms of the unit ball. This would allow for the exploration of diagrams for non-convex sets, such as the class of simply connected or star-shaped domains.
    \item \textbf{Further functionals}: Applying the method to more quantities, such as the eigenvalues of higher-order operators or the Cheeger constant, which are notoriously difficult to optimize numerically.
    \item \textbf{Other frameworks}: The underlying idea of sampling the image of a function can be transferred to other research areas, such as periodic homogenization theory, to systematically explore the set of achievable effective material properties induced by different microstructures.
\end{itemize}

{\bf Acknowledgments:} Eloi Martinet would like to thank Charles Dapogny and Leon Bungert for fruitful discussions.

This work was finalized during a stay of Ilias Ftouhi at the Isaac Newton Institute for Mathematical Sciences, Cambridge, in February 2026 during the program “Geometric spectral theory and applications”.  The author is grateful to the Institute for the very good and stimulating atmosphere. This work was partially supported by EPSRC grant EP/Z000580/1.

The authors acknowledge the use of ChatGPT (OpenAI) to assist in generating preliminary code templates. All generated outputs were thoroughly reviewed, validated, and, where necessary, revised by the authors. ChatGPT was also used for limited language refinement; however, it was not used to generate any substantive scientific content in this manuscript.

\bibliographystyle{abbrv}
\bibliography{biblio.bib}

\end{document}